\documentclass[12pt]{article}

\input{xypic.tex} \xyoption{all}
\usepackage{amsthm,amssymb,latexsym}
\usepackage{graphics}

\newcommand{\real}{\mathbb{R}} 
\newcommand{\Z}{\mathbb{Z}} 
\newcommand{\p}[1]{{\mathbb{P}^{#1}}} \newcommand{\op}[1]{{\cal O}_{\mathbb{P}^{#1}}}
\newcommand{\oh}{{\cal O}}  \newcommand{\cV}{{\cal V}}
 
\newcommand{\im}{{\rm Im}~} \newcommand{\pn}{{\mathbb{P}^n}}
\newcommand{\simto}{\stackrel{\sim}{\to}} \newcommand{\opn}{{\cal O}_{\mathbb{P}^n}}

\newtheorem{theorem}{Theorem}
\newtheorem{proposition}[theorem]{Proposition}
\newtheorem{lemma}[theorem]{Lemma}
\newtheorem{corollary}[theorem]{Corollary}
\newtheorem*{remark}{{\bf Remark}}
\newtheorem{example}{Example}
\newtheorem*{definition}{{\bf Definition}}
\newtheorem*{conjecture}{{\bf Conjecture}}

\begin{document}

\title{Instanton sheaves on complex projective spaces}
\author{Marcos Jardim \\ IMECC - UNICAMP \\
Departamento de Matem\'atica \\ Caixa Postal 6065 \\
13083-970 Campinas-SP, Brazil}

\maketitle

\begin{abstract}
We study a class of torsion-free sheaves on complex projective
spaces which generalize the much studied mathematical instanton
bundles. Instanton sheaves can be obtained as cohomologies of
linear monads and are shown to be semistable if its rank is not too large,
while semistable torsion-free sheaves satisfying certain cohomological
conditions are instanton. We also study a few examples of moduli spaces of 
instanton sheaves.

\end{abstract}

\tableofcontents

\newpage
\baselineskip18pt


\section*{Introduction} \label{intro}

The study of vector bundles and more general sheaves on complex projective
spaces has been a topic of great interest to algebraic geometers for many
years, see for instance the excellent book by Okonek, Schneider and
Spindler \cite{OSS} and Hartshorne's problem list \cite{Ha0}. In this paper
we concentrate on a particular class of sheaves defined as follows, generalizing
the concept of admissible sheaves on $\p3$ due to Manin \cite{Ma}, see also 
\cite{FJ2}. 

\begin{definition}
An {\em instanton sheaf} on $\pn$ ($n\geq2$) is a torsion-free
coherent sheaf $E$ on $\pn$ with $c_1(E)=0$ satisfying the following
cohomological conditions:
\begin{enumerate}
\item for $n\geq2$, $H^0(E(-1))=H^n(E(-n))=0$;
\item for $n\geq3$, $H^1(E(-2))=H^{n-1}(E(1-n))=0$;
\item for $n\geq4$, $H^p(E(k))=0$, $2\leq p\leq n-2$ and $\forall k$;
\end{enumerate}
The integer $c=-\chi(E(-1))$ is called the charge of $E$.
\end{definition}

If $E$ is a rank $2m$ locally-free instanton sheaf on $\mathbb{P}^{2m+1}$
of trivial splitting type (i.e. there exist a line $\ell\subset\mathbb{P}^{2m+1}$ such that
$E|_{\ell}\simeq\oh_{\ell}^{2m}$), then $E$ is a {\em mathematical instanton bundle}
as originally defined by Okonek and Spindler \cite{OS}. There is an extensive
literature on such objects, see for instance \cite{AO,SpT}. The nomenclature
is motivated by gauge theory: mathematical instanton bundles on $\mathbb{P}^{2n+1}$
correspond to quaternionic instantons on $\mathbb{P}\mathbb{H}^{n}$ \cite{Sa}.

The goal of this paper is to extend the discussion in two directions:
the inclusion of even-dimensional projective spaces and the analysis of
more general sheaves, allowing non-locally-free sheaves of arbitrary rank.
Such extension is motivated by the concept that in order to better understand
moduli spaces of stable vector bundles over a projective variety one must also consider
semistable torsion-free sheaves \cite{Ha0}. It turns out that many of the well-known
results regarding mathematical instanton bundles on $\mathbb{P}^{2n+1}$ generalize
in sometimes surprising ways to more general instanton sheaves.

The paper is organized as follows. We start by studying linear monads and their
cohomologies in Section \ref{s1}, spelling out criteria to decide whether the cohomology
of a given monad is torsion-free, reflexive or locally-free. We then show that 
every instanton sheaf is the cohomology of a linear monad, and that rank $r$
instanton sheaves on $\pn$ exist if and only if $r\geq n-1$. Further properties
of instanton sheaves are also studied in Section \ref{s2}.

The bulk of the paper lies in Section \ref{ss}, where we analyze the semistability
(in the sense of Mumford-Takemoto) of instanton sheaves. It is shown, for instance,
that every rank $r\le 2n-1$ locally-free instanton sheaf on $\pn$ is semistable, while
every rank $r\le n$ reflexive instanton sheaf on $\pn$ is semistable. We also determine
when a semistable torsion-free sheaf on $\pn$ is an instanton sheaf, showing for instance
that every semistable torsion-free sheaf on $\p2$ is an instanton sheaf.

In Section \ref{simple} it is shown that every rank $n-1$ instanton sheaf on $\pn$ is
simple, generalizing a result of Ancona and Ottaviani for mathematical instanton bundles
\cite[Proposition 2.11]{AO}. We then conclude in Section \ref{moduli} with a few results concerning the
moduli spaces of instanton sheaves.

It is also worth noting that Buchdahl has studied monads over arbitrary blow-ups
of $\p2$ \cite{B} while Costa and Mir\'o-Roig have initiated the study of
locally-free instanton sheaves over smooth quadric hypersurfaces within
$\pn$ \cite{CMR}, obtaining some results similar to ours. Many of the results here
obtained are also valid for instanton sheaves suitably defined over projective
varieties with cyclic Picard group, see \cite{JMR}.

\paragraph{Notation.}
We work over an algebraically closed field $\mathbb{F}$ of characteristic zero.
It might be interesting from the algebraic point of view to study how the results
here obtained generalize to finite fields. Throughout this paper, $U$, $V$ and $W$
are finite dimensional vector spaces over the fixed field $\mathbb{F}$, and we use
$[x_0:\cdots:x_n]$ to denote homogeneous coordinates on $\pn$. If $E$ is a sheaf
on $\pn$, then $E(k)=E\otimes\opn(k)$, as usual; by $H^p(E)$ we actually mean
$H^p(\pn,E)$ and $h^p(E)$ denotes the dimension of $H^p(\pn,E)$.

\bigskip

\paragraph{Acknowledgment.}
The author is partially supported by the FAEPEX grants number 1433/04 and 1652/04,
and the CNPq grant number 300991/2004-5. 
We thank Giorgio Ottaviani and Rosa Maria Mir\'o-Roig for their valuable
comments on the first version of this paper. This material was the topic
of a two-week course during the 2005 Summer Program at IMPA; we thank
Eduardo Esteves for the invitation and IMPA for the financial support.


\section{Monads} \label{s1}

Let $X$ be a smooth projective variety. A {\em monad} on $X$ is a
complex $V_{\bullet}$ of the following form:
\begin{equation} \label{m1}
\cV_{\bullet} ~ : ~ 0 \to V_{-1} \stackrel{\alpha}{\longrightarrow}
V_{0} \stackrel{\beta}{\longrightarrow} V_{1} \to 0
\end{equation}
which is exact on the first and last terms. Here, $V_k$ are
locally free sheaves on $X$. The sheaf $E=\ker\beta/\im\alpha$
is called the cohomology of the monad $\cV_{\bullet}$.

Monads were first introduced by Horrocks, who has shown that every
rank 2 locally free sheaf on $\p3$ can be obtained as the cohomology
of a monad where $V_k$ are sums of line bundles \cite{Ho}.

In this paper, we will focus on the so-called {\em linear monads}
on $\pn$, which are of the form:
\begin{equation} \label{spmo}
0\to  V\otimes\opn(-1) \stackrel{\alpha}{\longrightarrow}
W\otimes\opn \stackrel{\beta}{\longrightarrow} U\otimes\opn(1) \to 0 ~,
\end{equation}
where $\alpha\in{\rm Hom}(V,W)\otimes H^0(\opn(1))$ is injective (as a sheaf map)
and $\beta\in{\rm Hom}(W,U)\otimes H^0(\opn(1))$ is surjective. The {\em degeneration
locus} $\Sigma$ of the monad (\ref{spmo}) consists of the set of points $x\in\pn$
such that the localized map $\alpha_x\in{\rm Hom}(V,W)$ is not injective;
in other words $\Sigma = {\rm supp}~\{{\rm coker}\alpha^*\}$.

Linear monads have appeared in a wide variety of contexts within algebraic
geometry, like the construction of locally free sheaves on complex projective
spaces and the study of curves in $\p3$ and surfaces in $\p4$, see for instance
\cite{F} and the references therein. Our main motivation to study them comes
from gauge theory; as it is well-known, linear monads on $\p2$ and $\p3$ are
closely related to instantons on $\real^4$ \cite{D}.

The existence of linear monads on $\pn$ has been completely classified by Fl\o ystad
in \cite{F}; let $v=\dim V$, $w=\dim W$ and $u=\dim U$. 

\begin{theorem}\cite{F}\label{exist}
Let $n\geq1$. There exists a linear monad on $\pn$ as above if and only
if at least one of the following conditions hold:
\begin{enumerate}
\item[(i)] $w \geq 2u+ n - 1$ and $w \geq v+u$;
\item[(ii)] $w \geq v + u + n$.
\end{enumerate}
If the conditions hold, there exists a linear monad whose degeneration
locus is a codimension $w-v-u+1$ subvariety.
\end{theorem}

In particular, if $v,w,u$ satisfy condition (2) above, then the degeneration
locus is empty.

\begin{remark}\rm
A similar classification result for linear monads over $n$-dimensional quadric
hypersurfaces within $\mathbb{P}^{n+1}$ has been proved by Costa and Mir\'o-Roig
in \cite{CMR}, by adapting Fl\o ystad's technique. 
\end{remark}

\begin{definition}
A coherent sheaf on $\pn$ is said to be a linear sheaf if it can be
represented as the cohomology of a linear monad.
\end{definition}

The goal of this section is to study linear sheaves, with their
characterization in mind. First, notice that if $E$ is the cohomology
of (\ref{spmo}) then
$$ {\rm rank}(E) = w - v - u ~~ , ~~ c_1(E) = v - u ~~ {\rm and} $$
$$ {\rm c}(E)= \left(\frac{1}{1-H}\right)^v\left(\frac{1}{1+H}\right)^{u} ~~. $$

\begin{proposition}\label{coho}
If $E$ is a linear sheaf on $\pn$, then:
\begin{enumerate}
\item[(i)] for $n\geq2$, $H^0(E(k))=H^0(E^*(k))=0$, $\forall k\leq-1$;
\item[(ii)]for $n\geq3$, $H^1(E(k))=0$, $\forall k\leq-2$;
\item[(iii)] for $n\geq4$, $H^p(E(k))=0$, $2\leq p\leq n-2$ and $\forall k$;
\item[(iv)] for $n\geq3$, $H^{n-1}(E(k))=0$, $\forall k\geq-n+1$;
\item[(v)] for $n\geq2$, $H^n(E(k))=0$ for $k\geq-n$;
\item[(vi)] for $n\geq2$, ${\cal E}xt^1(E,\opn)={\rm coker}\alpha^*$ and
${\cal E}xt^p(E(k),\opn)=0$ for $p\geq2$ and all $k$.
\end{enumerate}\end{proposition}

In particular, note that linear sheaves have natural cohomology in
the range $-n\leq k\leq-1$, i.e. for the values of $k$ in this range at
most one of the cohomology groups $H^q(E(k))$ is nontrivial. It also follows
that $c=-\chi(E(-1))=h^1(E(-1))$.

Every rank $2n$ locally-free sheaf on $\mathbb{P}^{2n+1}$ with total Chern
class $c(E)=(1+H^2)^{-c}$ and natural cohomology in the rank $-2n-1\le k \le0$
is linear \cite{OS}, and therefore satisfy the stronger conclusion of the
Proposition above. However, not all sheaves on $\pn$ having natural cohomology
in the range $-n\leq k\leq-1$ are linear, the simplest example being
$\Omega^1_{\pn}(n)$.

\begin{proof}
Assume that $E$ is the cohomology of the monad (\ref{spmo}).
The kernel sheaf $K=\ker\beta$ is locally-free, and one has
the sequences $\forall k$:
\begin{equation} \label{ker1}
0 \to K(k) \to W\otimes\opn(k) \stackrel{\beta}{\longrightarrow}
U\otimes\opn(k+1) \to 0 ~~ {\rm and}
\end{equation}
\begin{equation} \label{ker2}
0 \to V\otimes\opn(k-1) \stackrel{\alpha}{\longrightarrow}
K(k) \to E(k) \to 0 ~~.
\end{equation}
From the first sequence, we see that: $H^p(K(k))=0$ for $p=0,1$ and $p+k\leq-1$; 
for $2\leq p\leq n-1$, $\forall k$; and for $p=n$ and $k\leq -n$.
The second sequence tell us that $H^p(K(k))\stackrel{\sim}{\to}H^p(E(k))$ for
$p=0$ and $k\leq0$, for $1\leq p\leq n-2$ and all $k$, and for
$p\geq n-1$ and $k\geq-n$. Putting these together, we obtain the
conditions $(i)$ through $(v)$ in the statement of the Proposition.

Dualizing sequence (\ref{ker1}), we obtain that $H^0(K^*(k))=0$
for $k\leq-1$. Now dualizing sequence (\ref{ker2}), we get, since
$K$ is locally-free:
\begin{equation} \label{ker3}
0 \to E^*(-k) \to K^*(-k) \stackrel{\alpha^*}{\to}
V\otimes\opn(-k+1) \to {\cal E}xt^1(E(k),\opn) \to 0 ~~ .
\end{equation}
Condition $(vi)$ and the second part of condition $(i)$ follow easily.
\end{proof}

Conversely, linear sheaves can be characterized by their cohomology:

\begin{theorem}\label{inst=mon}
If $E$ is a torsion-free sheaf on $\pn$ satisfying:
\begin{enumerate}
\item[(i)] for $n\geq2$, $H^0(E(-1))=H^n(E(-n))=0$;
\item[(ii)] for $n\geq3$, $H^1(E(-2))=H^{n-1}(E(1-n))=0$;
\item[(iii)] for $n\geq4$, $H^p(E(k))=0$, $2\leq p\leq n-2$ and $\forall k$;
\end{enumerate}
then $E$ is linear, and can be represented as the cohomology of the monad:
\begin{eqnarray} 
\label{m2} 0 & \to & H^1(E\otimes\Omega_{\pn}^2(1))\otimes\opn(-1) \to \\
\nonumber & \to & H^1(E\otimes\Omega^1_{\pn})\otimes\opn \to 
H^1(E(-1))\otimes\opn(1) \to 0 ~~ .
\end{eqnarray}\end{theorem}

\begin{proof}
Given a hyperplane $\wp\subset\pn$, consider the restriction sequence:
$$ 0\to E(k-1) \to E(k) \to E(k)|_{\wp} \to 0 ~~.$$
Clearly, $H^0(E(-1))=0$ implies that $H^0(E(k))=0$ for $k\le-1$, while
$H^n(E(-n))=0$ forces $H^n(E(k))=0$ for $k\ge-n$.

Since $H^0(E(-1))=H^1(E(-2))=0$, it follows that $H^0(E(-1)|_{\wp})=0$,
hence $H^0(E(k)|_{\wp})=0$ for $k\le-1$. So we have the sequence:
$$ 0 \to H^1(E(k-1)) \to H^1(E(k)) ~~,~~ {\rm for} ~~ k\le-2 $$
thus by induction $H^1(E(k))=0$ for $k\le-2$.

Since $H^n(E(-n))=H^{n-1}(E(1-n))=0$, it follows that $H^{n-1}(E(1-n)|_{\wp})=0$,
hence by further restriction $H^{n-1}(E(k)|_{\wp})=0$ for $k\ge1-n$. So we
have the sequence:
$$ H^{n-1}(E(k-1)) \to H^{n-1}(E(k)) \to 0  ~~,~~ {\rm for} ~~ k\ge1-n $$
thus by induction $H^{n-1}(E(k))=0$ for $k\ge1-n$.

The key ingredient of the proof is the Beilinson spectral sequence \cite{OSS}:
for any coherent sheaf $E$ on $\pn$ there exists a spectral sequence $\{E^{p,q}_r\}$
whose $E_1$-term is given by ($q=0,\dots,n$ and $p=0,-1,\dots,-n$):
$$ E_1^{p,q} = H^q(E\otimes\Omega_{\pn}^{-p}(-p))\otimes \opn(p) $$
which converges to
$$ E^i = \left\{ \begin{array}{l} 
E ~,~ {\rm if} ~ p+q=0 \\ 0 ~ {\rm otherwise} \end{array} \right. ~~ . $$

Applying the Beilinson spectral sequence to $E(-1)$, we must show that
\begin{equation}\label{vhq}
H^q(E(-1)\otimes\Omega_{\pn}^{-p}(-p))=0 ~~ {\rm for}~q\neq1 ~~
{\rm and~for}~ q=1,~p\leq-3 ~ .
\end{equation}
It then follows that the Beilinson spectral sequence degenerates at the $E_2$-term
and the monad
\begin{eqnarray} 
\label{m3} 0 & \to & H^1(E(-1)\otimes\Omega_{\pn}^2(2))\otimes\opn(-2) \to \\
\nonumber & \to & H^1(E(-1)\otimes\Omega^1_{\pn}(1))\otimes\opn(-1) \to 
H^1(E(-1))\otimes\opn \to 0
\end{eqnarray}
has $E(-1)$ as its cohomology. Tensoring (\ref{m3}) by $\opn(1)$, we conclude
that $E$ is the cohomology of (\ref{m2}), as desired.

The claim (\ref{vhq}) follows from repeated use of the exact sequence
$$ H^q(E(k))^{\oplus m} \to H^q(E(k+1)\otimes\Omega_{\pn}^{-p-1}(-p-1)) \to $$
\begin{equation} \label{e(k)euler-hom}
\to H^{q+1}(E(k)\otimes\Omega_{\pn}^{-p}(-p)) \to H^{q+1}(E(k))^{\oplus m}
\end{equation}
associated with Euler sequence for $p$-forms on $\pn$ twisted by $E(k)$:
\begin{equation} \label{e(k)euler}
0 \to E(k)\otimes\Omega_{\pn}^{-p}(-p) \to E(k)^{\oplus m} \to
E(k)\otimes\Omega_{\pn}^{-p-1}(-p) \to 0 ~,
\end{equation}
where $q=0,\dots,n$ , $p=-1,\dots,-n$ and
$m=\left(\begin{array}{c}n+1\\-p\end{array}\right)$.

For instance, it is easy to see  that:
$$ H^0(E(k)\otimes\Omega_{\pn}^{-p}(-p))=0 ~~ {\rm for~all} ~ p ~ {\rm and} ~ k\leq-1 ~ ; $$
$$ H^q(E(-1)\otimes\Omega_{\pn}^{n}(n))=H^q(E(-2))=0 ~~ {\rm for~all} ~ q ~ ; $$
$$ H^q(E(-1))=0 ~~ {\rm for~all} ~ q\neq1 ~ ; $$
$$ H^n(E(k)\otimes\Omega_{\pn}^{-p}(-p))=0 ~~ {\rm for~all} ~ p ~ {\rm and} ~ k\geq-n ~. $$
Setting $q=n-1$, we also obtain:
$$ H^{n-1}(E(k)\otimes\Omega_{\pn}^{-p}(-p))=0 ~~ {\rm for} ~ p\geq-n+1
~ {\rm and} ~ k\geq-n-1 ~, $$
and so on.
\end{proof}

Clearly, the cohomology of a linear monad is always coherent, but more can be
said if the codimension of the degeneration locus of $\alpha$ is known.

\begin{proposition}
Let $E$ be a linear sheaf. 
\begin{enumerate}
\item[(i)] $E$ is locally-free if and only if its degeneration locus is empty;
\item[(ii)] $E$ is reflexive if and only if its degeneration locus is a subvariety
of codimension at least 3;
\item[(iii)] $E$ is torsion-free if and only if its degeneration locus is a subvariety
of codimension at least 2.
\end{enumerate} \end{proposition}

\begin{proof}
Let $\Sigma$ be the degeneration locus of the linear sheaf $E$. From Proposition
\ref{coho}, we know that ${\cal E}xt^p(E,\opn)=0$ for $p\geq2$ and
$$ \Sigma = {\rm supp}~{\cal E}xt^1(E,\opn) = 
\{ x\in\pn ~ | ~ \alpha_x ~ {\rm is~not~injective} ~ \} ~~ . $$

The first statement is clear; so it is now enough to argue that
$E$ is torsion-free if and only if $\Sigma$ has codimension at least
2 and that $E$ is reflexive if and only if $\Sigma$ has codimension
at least 3.

Recall that the $m^{\rm th}$-singularity set of a coherent sheaf $\cal F$
on $\pn$ is given by:
$$ S_m({\cal F}) = \{ x\in\pn ~|~ dh({\cal F}_x) \geq n-m \} $$
where $dh({\cal F}_x)$ stands for the homological dimension of
${\cal F}_x$ as an ${\cal O}_x$-module:
$$ dh({\cal F}_x) = d ~~~ \Longleftrightarrow ~~~
\left\{ \begin{array}{l}
{\rm Ext}^d_{{\cal O}_x}({\cal F}_x,{\cal O}_x) \neq 0 \\
{\rm Ext}^p_{{\cal O}_x}({\cal F}_x,{\cal O}_x) = 0 ~~ \forall p>d
\end{array} \right. $$

In the case at hand, we have that $dh(E_x) = 1$ if $x\in\Sigma$,
and $dh(E_x) = 0$ if $x\notin\Sigma$. Therefore
$S_0(E)=\cdots=S_{n-2}(E)=\emptyset$, while $S_{n-1}(E)=\Sigma$.
It follows that \cite[Proposition 1.20]{ST}:
\begin{itemize}
\item[(i)] if ${\rm codim}~\Sigma \geq 2$, then $\dim S_m(E)\leq m-1$ for all $m<n$,
hence $E$ is a locally 1$^{\rm st}$-syzygy sheaf;
\item[(ii)]  if ${\rm codim}~\Sigma \geq 3$, then $\dim S_m(E)\leq m-2$ for all $m<n$,
hence $E$ is a locally 2$^{\rm nd}$-syzygy sheaf.
\end{itemize}
The desired statements follow from the observation that $E$ is
torsion-free if and only if it is a locally 1$^{\rm st}$-syzygy sheaf, 
while $E$ is reflexive if and only if it is a locally
2$^{\rm nd}$-syzygy sheaf \cite[p. 148-149]{OSS}. 
\end{proof}


\paragraph{A splitting criterion for locally-free linear sheaves.}
Given a coherent sheaf $E$ on $\pn$, we define
$$ H^p_*(E) = \bigoplus_{k\in\Z} H^p(E(k)) $$
which has the structure of a graded module
over $S^n=\bigoplus_{k\in\Z} H^p(\oh(k))$. Kumar, Peterson and
Rao prove the following result \cite{KPR}:

\begin{theorem}\label{mkr}
Let $E$ be a rank $r$ locally-free sheaf on $\pn$, $n\geq4$.
\begin{enumerate}
\item[(i)] If $n$ is even and $r\leq n-1$, then $E$ splits as a sum of line bundles
if and only if $H^p_*(E)=0$ for $2\leq p\leq n-2$.
\item[(ii)] If $n$ is odd and $r\leq n-2$, then $E$ splits as a sum of line bundles
if and only if $H^p_*(E)=0$ for $2\leq p\leq n-2$.
\end{enumerate}\end{theorem}

Thus we obtain as an easy consequence of (3) in Proposition
\ref{coho} and the previous theorem:

\begin{corollary}
Let $E$ be a rank $r$ locally-free linear sheaf on $\pn$.
\begin{enumerate}
\item[(i)] If $n$ is even and $r\leq n-1$, then $E$ splits as a sum of line bundles.
\item[(ii)] If $n$ is odd and $r\leq n-2$, then $E$ splits as a sum of line bundles.
\end{enumerate}\end{corollary}

This means that linear monads are not useful to produce locally free sheaves
of low rank on $\pn$, one of the problems suggested by Hartshorne in \cite{Ha0}
and still a challenge in the subject.

Let us also point out that Kumar, Peterson and Rao's result is optimal, in the
sense that there exist rank $2n$ locally-free sheaves on $\mathbb{P}^{2n+1}$
and on $\mathbb{P}^{2n}$ ($n\geq2$) satisfying $H^p_*(E)=0$ for $2\leq p\leq 2n-1$
which are stable, and hence do not split as a sum of line bundles. Moreover, as
we will see in examples below, it does not generalize to reflexive or torsion-free
sheaves either.

Thus linear monads will only produced interesting locally-free sheaves when $r=w-v-u\geq n$
if $n$ is even and when $r=w-v-u\geq n-1$ if $n$ is odd. However, we can still expect
to use monads to produce interesting {\em reflexive} sheaves of low rank,
see Example \ref{ex4} below.


\section{Basic properties of instanton sheaves}\label{s2}

It follows from Proposition \ref{coho} that the cohomology of a linear monad
with $v=u$ and $w-2v\geq 1$ is a torsion-free instanton sheaf of rank $w-2v$
and charge $v$. In particular, by Fl\o ystad's theorem, there are rank $r$
instanton sheaves on $\pn$ for each $r\ge n-1$. Moreover, for every $r\geq2n$,
there are rank $r$ {\em locally-free} instanton sheaves on $\mathbb{P}^{2n+1}$
or $\mathbb{P}^{2n}$.

Moreover, by virtue of Theorem \ref{inst=mon}, every rank $r$ instanton sheaf
of charge $c$ is the cohomology of a monad of the form:
\begin{equation}\label{instmonad}
0 \to \opn(-1)^{\oplus c} \stackrel{\alpha}{\to} \opn^{\oplus r+2c}
\stackrel{\beta}{\to} \opn(1)^{\oplus c} \to 0
\end{equation}
for some injective map $\alpha$ degenerating in codimension at least 2 and some
surjective map $\beta$, such that $\beta\alpha=0$. It follows easily from Fl\o ystad's
theorem that there are no instanton sheaves on $\pn$ of rank $r\leq n-2$.

\begin{corollary}\label{edual}
If $E$ is an instanton sheaf then:
\begin{itemize}
\item[(i)] $H^0(E^*(k))=0$, $\forall k\leq-1$;
\item[(ii)] ${\cal E}xt^p(E,\opn)=0$ for $p\geq2$;
\item[(iii)] ${\rm Ext}^p(E,E)=0$ for $p\geq3$.
\end{itemize} \end{corollary}
\begin{proof}
The first two statements follow easily from Proposition \ref{coho} and the fact that
every instanton sheaf is linear. For the last statement, let $K=\ker\beta$; taking
the sequence:
$$ 0\to \opn(-1)^{\oplus c} \to K \to E \to 0 ~~, $$
we obtain:
\begin{equation}\label{extsqc1}
{\rm Ext}^p(E,\opn(-1)^{\oplus c}) \to {\rm Ext}^p(E,K) \to 
{\rm Ext}^p(E,E) \to {\rm Ext}^{p+1}(E,\opn(-1)^{\oplus c})
\end{equation}
Thus ${\rm Ext}^p(E,K)\simeq{\rm Ext}^p(E,E)$ for all $p\geq2$.
Now from the sequence
$$ 0 \to K \to \opn^{\oplus r+2c} \to \opn(1)^{\oplus c} \to 0 $$
we obtain:
\begin{equation}\label{extsqc2}
{\rm Ext}^p(E,\opn(1)^{\oplus c})\to {\rm Ext}^{p+1}(E,K) \to
{\rm Ext}^{p+1}(E,\opn^{\oplus r+2c}) ~~.
\end{equation}
Hence ${\rm Ext}^{p}(E,K)=0$ for $p\geq3$, and the result follows.
\end{proof}

Furthermore, ${\rm Ext}^1(E,E)$ is completely determined by the map
$\alpha$, while ${\rm Ext}^2(E,E)$ is completely determined by the map
$\beta$. Indeed, first note that the maps $\alpha$ and $\beta$ induce
linear maps:
$$ {\rm Ext}^1\alpha:{\rm Ext}^1(E,\opn(-1)^{\oplus c}) \to {\rm Ext}^1(E,K) $$
$$ {\rm Ext}^1\beta:{\rm Ext}^1(E,\opn^{\oplus r+2c}) \to {\rm Ext}^1(E,\opn(1)^{\oplus c}) $$
Setting $p=1$ on (\ref{extsqc1}) and (\ref{extsqc2}) we get:
$$ {\rm Ext}^1(E,E) = {\rm coker}~\{{\rm Ext}^1\alpha\} ~~{\rm and}~~
{\rm Ext}^2(E,E) = {\rm coker}~\{{\rm Ext}^1\beta\} $$
In particular, note that ${\rm Ext}^2(E,E)=0$ is an open condition on \linebreak
${\rm Hom}(\opn^{\oplus r+2c},\opn(1)^{\oplus c})$.

Given two linear sheaves, one can produce a new instanton sheaf of
higher rank using the following result:

\begin{proposition}\label{ext}
An extention $E$ of linear sheaves $F'$ and $F''$
$$ 0 \to F' \to E \to F'' \to 0 $$
is also a linear sheaf. Moreover, if $c_1(F')=-c_1(F'')$, then
$E$ is instanton.
\end{proposition}

\begin{proof}
The desired statement follows easily from the associated 
sequences of cohomology:
$$ H^q(F'(k)) \to H^q(E(k)) \to H^q(F''(k)) ~ , ~ \forall q=0,\dots,n ~ ,$$
so that $H^q(E(k))$ vanishes whenever $H^q(F'(k))$ and $H^q(F''(k))$ do.
Note that $E$ is classified by ${\rm Ext}^1(F'',F')$.
\end{proof}

\begin{proposition}
If $E$ is a locally-free instanton sheaf on $\pn$,
then $E^*$ is also instanton.
\end{proposition}

\begin{proof}
The statement is an easy consequence of Serre duality. In fact, if $E$ arises
as the cohomology of the monad
$$ 0\to  V\otimes\opn(-1) \stackrel{\alpha}{\longrightarrow}
W\otimes\opn \stackrel{\beta}{\longrightarrow} U\otimes\opn(1) \to 0 ~, $$
then $E^*$ is the cohomology of the dual monad
$$ 0\to  U^*\otimes\opn(-1) \stackrel{\beta^*}{\longrightarrow}
W^*\otimes\opn \stackrel{\alpha^*}{\longrightarrow} V^*\otimes\opn(1) \to 0 ~. $$
\end{proof}

In general, if $E$ is not locally-free, its dual might not be an
instanton sheaf, see example \ref{dnoi} below. However, the dual of
every semistable sheaf on $\p2$ is instanton.

\begin{proposition}
Let $\wp$ be a hyperplane in $\pn$. The restriction $E|_{\wp}$ of an
instanton sheaf $E$ on $\pn$ is also an instanton sheaf, and the restriction
map $\rho:H^1(E(-1))\to H^1(E(-1)|_{\wp})$ is an isomorphism.
\end{proposition}
\begin{proof}
Follows easily from the definition and the exact sequence
$$ 0\to E(k-1) \stackrel{\sigma}{\to} E(k) \to E|_{\wp}(k) \to 0 ~~.$$
\end{proof}

\paragraph{Instanton sheaves and mathematical instanton bundles.}
Mathematical instanton bundles have been defined in \cite{OS} as
a rank $2m$ locally-free sheaf on $\mathbb{P}^{2m+1}$ satisfying
the following conditions:
\begin{itemize}
\item $c(E) = \left(\frac{1}{1-H^2}\right)^c = (1+H^2+H^4+\cdots)^c$;
\item $E$ has natural cohomology in the range $-2m-1\leq k\leq0$;
\item $E$ is simple;
\item $E$ has trivial splitting type (i.e. there exist a line
$\ell\subset\mathbb{P}^{2m+1}$ such that $E|_{\ell}\simeq\oh_{\ell}^{2m}$).
\end{itemize} 
It was later shown by Ancona and Ottaviani have shown that the
simplicity assumption is redundant \cite[Proposition 2.11]{AO}:
every rank $2m$ locally-free sheaf on $\mathbb{P}^{2m+1}$ satisfying
the first two conditions is simple (in fact, more is true, see Lemma
\ref{simp} below). The last condition is also redundant for $m=1$
and for $k=1$; there are however rank $2m$ locally-free sheaf on
$\mathbb{P}^{2m+1}$ satisfying the first two conditions which are
not of trivial splitting type.

Every mathematical instanton bundle as above can be represented as the
cohomology of a linear monad with $v=u=c$ and $w=2m+2c$ \cite{OS}, hence
it is a rank $2m$ locally-free instanton sheaf on $\mathbb{P}^{2m+1}$
of charge $c$.

Conversely, a rank $2m$ locally-free instanton sheaf $E$ on $\mathbb{P}^{2m+1}$
is a mathematical instanton bundle as above if $H^0(E)=0$ and it is of trivial
splitting type. We show that the vanishing of $H^0(E)$ is automatic.

\begin{proposition}\label{h0}
If $E$ is a rank $n-1$ instanton sheaf on $\pn$, then
$H^0(E)=0$. If $E$ is locally-free, then $H^0(E^*)=0$.
\end{proposition}
\begin{proof}
Let $E$ be a rank $n-1$ locally-free instanton sheaf on $\pn$,
and assume that $H^0(E)\neq0$. So let $Q=E/\opn$ and
consider the sequence:
$$ 0 \to \opn \to E \to Q \to 0 ~~.$$
Note that $c_1(Q)=0$, $H^0(Q(k))=0$ for $k\leq-1$, $H^n(Q(k))=0$
for $k\geq-n$ and $H^q(Q(k))=H^q(E(k))$ for $1\leq q\leq n-1$. It
follows that $Q$ must be a rank $n-2$ instanton sheaf, which
cannot exist by Fl\o ystad's theorem.

If $E$ is locally-free, then $E^*$ is also instanton and
the vanishing of $H^0(E^*)$ follows by the same argument.
\end{proof}

In particular, there are no rank $2m-1$ locally-free instanton sheaves
on $\mathbb{P}^{2m}$. Indeed, by Theorem \ref{mkr}, any such sheaf must
split as a sum of line bundles, and this contradicts $c_1(E)=0$ and
$H^0(E)=0$.

We also point out that there are rank $n-1$ properly reflexive instanton sheaves
$E$ on $\pn$ for which $H^0(E^*)\neq0$, see Example \ref{p5} below.


\section{Semistability of instanton sheaves}\label{ss}

Recall that a torsion-free sheaf $E$ on $\pn$ is said to be
{\em semistable} if for every coherent subsheaf
$0\neq F\hookrightarrow E$ we have
$$ \mu(F)=\frac{c_1(F)}{{\rm rk}(F)} \leq \frac{c_1(E)}{{\rm rk}(E)}=\mu(E) ~ . $$
Furthermore, if for every coherent subsheaf $0\neq F\hookrightarrow E$
with $0<{\rm rk}(F)<{\rm rk}(E)$ we have
$$ \frac{c_1(F)}{{\rm rk}(F)} < \frac{c_1(E)}{{\rm rk}(E)} ~ , $$
then $E$ is said to be {\em stable}.
A sheaf is said to be properly semistable if it is semistable but not
stable. It is also important to remember that $E$ is (semi)stable if
and only if $E^*$ and $E(k)$ are.

For any given torsion-free sheaf $E$ of rank $r$, there is an uniquely
determined integer $k_E$ such that
$$ c_1(E(k_E)) = c_1(E) + rk_E \in \{0,-1,\cdots,-r+1\} ~~ ; $$
$E_\eta=E(k_E)$ is called the normalization of $E$. A sheaf $E$ is said to be
normalized if $-r+1\leq c_1(E)\leq0$.

\begin{lemma} ({\rm \cite[p. 167]{OSS}}) \label{r2}
Let $E$ be a normalized torsion-free sheaf on $\pn$.
If $E$ is stable then  $H^0(E)=0$ and
\begin{itemize}
\item $H^0(E^*)=0$ if $c_1(E)=0$;
\item $H^0(E^*(-1))=0$ if $c_1(E)<0$.
\end{itemize}
If $E$ is semistable then $H^0(E(-1))=H^0(E^*(-1))=0$.
\end{lemma}

For sheaves of rank 2 or 3, the above necessary criteria turns out
to be also sufficient, as we recall in the next two lemmas. We fix $n\geq2$.

\begin{lemma} \label{rk2tf}
Let $E$ be a normalized rank 2 torsion-free sheaf $E$ on $\pn$.
If $c_1(E)=0$, then:
\begin{itemize}
\item $E$ is stable if and only if $H^0(E)=H^0(E^*)=0$;
\item $E$ is semistable if and only if $H^0(E(-1))=H^0(E^*(-1))=0$.
\end{itemize}
If $c_1(E)=-1$, then $E$ is stable if and only if it is semistable if and
only if $H^0(E)=H^0(E^*(-1))=0$.
\end{lemma}
\begin{proof}
For $E$ being reflexive, this result is in \cite[p. 166]{OSS}. In general,
simply note that $E^*$ is reflexive and use the result just mentioned.
\end{proof}

As an easy consequence of Proposition \ref{coho} and Lemma \ref{rk2tf} we have:

\begin{proposition} \label{m+rk2+rk3}
If $E$ is a rank 2 instanton sheaf on $\pn$ ($n=2,3$), then $E$ is semistable.
It is stable if $H^0(E)=0$.
\end{proposition}

This result does not generalizes for higher rank, even we restrict ourselves to
locally-free instanton sheaves, as we show in the example below.

\begin{example}\label{unstable.p2}\rm
For each $n\geq 2$, there are rank $2n$ locally-free instanton sheaves on $\pn$
which are not semistable. Indeed, by Fl\o ystad's theorem, there is a linear monad:
\begin{equation} \label{def-F}
0 \to \opn(-1)^{\oplus a+1} \stackrel{\alpha}{\to} \opn^{\oplus n+2a+1}
\stackrel{\beta}{\to} \opn(1)^{\oplus a} \to 0 ~~ (a\ge1) ~~,
\end{equation}
whose cohomology $F$ is a locally-free sheaf of rank $n$ on $\pn$ and $c_1(F)=1$.
The dual $F^*$ is a locally-free sheaf of rank $n$ on $\pn$ and $c_1(F^*)=-1$. Any
extention of $E$ of $F^*$ by $F$:
$$ 0\to F\to E\to F^*\to 0 $$
is a rank $2n$ locally-free instanton sheaf which is clearly not semistable. Furthermore,
one can adjust the parameter $a$ depending on $n$ to ensure the existence of nontrivial
extensions.
\end{example}

For instanton sheaves of higher rank, the best statement one can have is the following:

\begin{theorem}\label{highrank}
Let $E$ be a rank $r$ instanton sheaf on $\pn$.
\begin{itemize}
\item If $E$ is reflexive and $r\le n+1$, then $E$ is semistable;
\item if $E$ is locally-free and $r\le 2n-1$, then $E$ is semistable.
\end{itemize}\end{theorem}

\begin{example}\rm
Note that the upper bound in the rank given in the second part of Theorem
\ref{highrank} is sharp, as seen in Example \ref{unstable.p2}. We now show
that there are rank $n+2$ reflexive instanton sheaves which are not semistable.

Indeed, let $X=\pn$, $n\ge3$. By Fl\o ystad's theorem \cite{F}, there is a
linear monad:
$$ 0 \to \opn(-1)^{\oplus n-2} \stackrel{\alpha}{\to} \opn^{\oplus n+1}
\stackrel{\beta}{\to} \opn(1) \to 0 $$
whose cohomology $F$ is a rank $2$ reflexive linear sheaf on $\pn$ and
$c_1(F)=n-3$.

Next, consider the rank $n$ locally free linear sheaf $G$
associated to the linear monad:
$$ 0 \to \opn(-1)^{\oplus a} \stackrel{\alpha}{\to} \opn^{\oplus 2n+2a-3}
\stackrel{\beta}{\to} \opn(1)^{\oplus n+a-3} \to 0 ~~ (a\ge 1). $$
Note that $c_1(G)=3-n$.

As in the previous example, an extension of $G$ by $F$ is a rank $n+2$ reflexive
instanton sheaf which is not semistable. The choice of a suitable value of the
parameter $a$ guarantees the existence of non-trivial extensions.
\end{example}

The proof of Theorem \ref{highrank} is based on Hoppe's criterion
\cite{Hop}: if $E$ is a rank $r$ reflexive sheaf on $\pn$ with $c_1(E)=0$
satisfying
$$ H^0(\Lambda^qE\!(-1))=0 ~~ for ~~ 1\leq q\leq r-1 $$
then $E$ is semistable. Indeed, assume $E$ is not semistable,
and let $F$ be a rank $q$ destabilizing sheaf with $c_1(F)=d>0$.
Then $\Lambda^qF=\opn(d)$, and the induced map
$\Lambda^qF\to\Lambda^qE$ yields a section in $H^0(\Lambda^qE\!(-d))$,
which forces $h^0(\Lambda^qE\!(-1))\neq0$. Similarly, it is also easy
to see that if
$$ H^0(\Lambda^qE)=0 ~~ for ~~ 1\leq q\leq r-1 $$
then $E$ is stable.

\noindent{\em Proof of Theorem \ref{highrank}.}
Every rank $r$ reflexive instanton sheaf on $\pn$ can be represented
as the cohomology of the monad (\ref{instmonad}).  Taking the sequence
$$ 0\to K\to \opn^{\oplus(r+2c)} \stackrel{\beta}{\to} \opn(1)^{\oplus c}\to0 ~~,$$
we consider the associated long exact sequence of exterior powers,
twisted by $\opn(-1)$:
$$ 0\to \Lambda^qK\!(-1)\to \Lambda^q(\opn^{\oplus(r+2c)})(-1) \to \cdots ~~.$$
Hence $H^0(\Lambda^qK\!(-1))=0$ for $1\leq q\leq r+c-1$. Now take the
sequence:
$$ 0\to \opn(-1)^{\oplus c}\to K\to E\to 0 ~~,$$
and consider the associated long exact sequence of symmetric powers,
twisted by $\opn(-1)$:
$$ 0\to \opn(-q-1)^{{c+q-1\choose q}} \to K(-q)^{{c+q-2\choose q-1} }\to
\cdots $$
$$ \to  \Lambda^{q-1}K\otimes\opn(-2)^{\oplus c} \to \Lambda^qK(-1)
\to \Lambda^qE(-1)\to 0 ~~.$$
Cutting into short exact sequences and passing to cohomology, we have
obtain that every reflexive instanton sheaf satisfies:
\begin{equation}\label{um}
H^0(\Lambda^pE(-1))=0~~ {\rm for} ~~ 1\leq p\leq n-1 ~~.
\end{equation}

It follows from (\ref{um}) that every rank $r\le n$ reflexive instanton
sheaf is semistable. If $E$ is a rank $n+1$ reflexive instanton sheaf, then
because $c_1(E)=0$:
$$ H^0(\Lambda^nE(-1))=H^0(E^*(-1))=0 ~~, $$
thus $E$ is also semistable.

Now if $E$ is locally-free,  the dual $E^*$ is also an instanton sheaf on $X$, so
\begin{equation}\label{dois}
H^0(\Lambda^q(E^*)(-1))=0 ~~ {\rm for} ~~ 1\leq q\leq n-1 ~~.
\end{equation}
But $\Lambda^p(E^*)\simeq\Lambda^{r-p}(E^*)$, since $\det(E)=\opn$; it follows
that:
\begin{eqnarray}
\nonumber H^0(\Lambda^pE(-1))=H^0(\Lambda^{r-p}(E^*)(-1))=0 &~~ {\rm for} ~~& 
1\leq r-p \leq n-1 \\
\label{tres}&\Longrightarrow& r-n+1\leq p \leq r-1 ~~.
\end{eqnarray}

Together, (\ref{dois}) and (\ref{tres}) imply that if $E$ is a rank $r\leq 2n-1$
locally-free instanton sheaf, then:
$$ H^0(\Lambda^pE(-1))=0~~ {\rm for} ~~ 1\leq p\leq 2n-2 $$
hence $E$ is semistable by Hoppe's criterion.
\hfill$\Box$

\begin{example}\rm
A similar result for the semistability of torsion-free instanton sheaves
beyond rank $2$ is unclear. However, it is easy to construct rank $n+1$
torsion-free instanton sheaves which are not semistable. Indeed, the
cohomology of the monad:
$$ 0 \to \opn(-1)^{\oplus n-1} \to \opn^{\oplus n+1} \to \opn(1) \to 0 $$
is of the form ${\cal I}_M(n-2)$, where ${\cal I}_M$ is the ideal sheaf
of a codimension 2 subvariety $M\hookrightarrow\pn$ \cite{F}.

On the other hand, there is a rank $n$ locally-free linear sheaf $F$ on
$\pn$ with $c_1(F)=2-n$ given by the cohomology of the monad:
$$ 0 \to \opn(-1)^{\oplus c} \to \opn^{\oplus 2c+2n-2} \to
\opn(1)^{\oplus c+n-2} \to 0 ~~. $$
Thus the sheaf $E$ given by the extention:
$$ 0 \to {\cal I}_M(n-2) \to E \to F \to 0 $$
is a rank $n+1$ torsion-free instanton sheaf which is not semistable.

In other words, Proposition \ref{m+rk2+rk3} is sharp on $\p2$, and
the reasonable conjecture seems to be that every rank $r=n-1,n$
torsion-free instanton sheaves on $\pn$ are semistable.
\end{example}

On the other hand, we have:

\begin{proposition} \label{nost}
For $r>(n-1)c$, there are no stable rank $r$ instanton sheaves on $\pn$ of
charge $c$.
\end{proposition}

In other words, every stable rank $r$ instanton sheaf on $\pn$ must be
of charge $c\geq r/(n-1)$, and there are properly semistable rank $r$
instanton sheaves for $n\leq r\leq 2n-1$.

\begin{proof}
Any rank $r$ instanton sheaf of charge $c$ is the cohomology of a monad of
the form:
$$ 0 \to \opn(-1)^{\oplus c} \stackrel{\alpha}{\to} \opn^{\oplus r+2c}
\stackrel{\beta}{\to} \opn(1)^{\oplus c} \to 0 $$
and note that
$$ H^0(E) \simeq H^0(\ker\beta) \simeq
\ker\{~H^0\beta~:~H^0(\opn^{\oplus r+2c})\to H^0(\opn(1)^{\oplus c})~\} $$
thus if $r>(n-1)c$, then $h^0(\opn^{\oplus r+2c})>h^0(\opn(1)^{\oplus c})$
and $H^0(E)\neq0$, so that $E$ is not stable by Lemma \ref{r2}.
\end{proof}

Let us now analyze the inverse question: are all semistable sheaves of degree
zero on $\pn$ instanton? The answer is positive for $n=2$, but there are
cohomological restrictions for $n\geq3$.

\begin{theorem} \label{s-s=>inst}
Let $E$ be a torsion-free sheaf on $\p2$ with $c_1(E)=0$.
If $E$ is semistable, then $E$ is instanton.
\end{theorem}
\begin{proof}
The semistability of $E$ and $E^*$ immediately implies that
$H^0(E(k))=H^0(E^*(k))=0$ for $k\leq-1$. If $E$ is a
locally-free sheaf, then via Serre duality $H^2(E(k))=0$
for $k\geq-2$, thus $E$ is instanton.

Now if $E$ is properly torsion-free, we consider the sequence:
\begin{equation}\label{tf}
0 \to E \to E^{**} \to Q \to 0
\end{equation}
where $Q=E^{**}/E$ is supported on a zero dimensional subscheme. Clearly,
$E^{**}$ is a semistable locally-free sheaf with $c_1(E)=0$, so it is
instanton by the previous paragraph. It follows from (\ref{tf}) that:
$$ H^0(E(k)) \hookrightarrow H^0(E^{**}(k))=0 ~~
{\rm for}~ k\leq-1 ~,~ {\rm and} $$
$$ H^2(E(k)) \stackrel{\simeq}{\rightarrow} H^2(E^{**}(k))=0 ~~
{\rm for}~ k\geq-2 ~, $$
so $E$ is also instanton.
\end{proof}

For $n\geq 3$, we have:

\begin{proposition} \label{s-s=>inst2}
If $E$ is a semistable locally-free sheaf on $\pn$ with $c_1(E)=0$
such that $H^1(E(-2))=H^{n-1}(E(1-n))=0$ and, for $n\geq4$, $H^p_*(E)=0$
for $2\leq p\leq n-2$, then $E$ is instanton.
\end{proposition}
\begin{proof}
If $E$ is semistable, then $H^0(E(k))=H^0(E^*(k))=0$ for $k\leq-1$,
hence $H^n(E(k))=0$ for $k\geq-n$ by Serre duality.
\end{proof}

As simple consequence of Proposition \ref{m+rk2+rk3} and Proposition
\ref{s-s=>inst2} we have:
\begin{itemize}
\item A rank 2 torsion-free sheaf on $\p3$ with $c_1(E)=0$ is instanton if and
only if it is semistable and $H^1(E(-2))=H^2(E(-2))=0$.
\item A rank 3 reflexive sheaf on $\p3$ with $c_1(E)=0$ is instanton if and
only if it is semistable and $H^1(E(-2))=H^2(E(-2))=0$.
\item A rank 3 reflexive sheaf on $\p4$ with $c_1(E)=0$ is instanton if and
only if it is semistable and $H^1(E(-2))=H^2_*(E)=H^3(E(-3))=0$.
\item A rank $4\le r \le 2n-1$ locally-free sheaf on $\pn$ ($n\geq3$) with $c_1(E)=0$
is instanton if and only if it is semistable and $H^1(E(-2))=H^{n-1}(E(1-n))=0$
and, for $n\geq4$, $H^p_*(E)=0$ for $2\leq p\leq n-2$.
\end{itemize}

\begin{remark}\rm
Since every Gieseker semistable torsion-free sheaf on $\pn$
is semistable \cite[p. 174]{OSS}, one can use the results above to decide
when a Gieseker semistable torsion-free sheaf on $\pn$ is instanton. It is
easy to see, however, that not all instanton sheaves are Gieseker semistable; indeed
if $E$ is an instanton sheaf satisfying $H^0(E)\neq0$, then $E$ is not
Gieseker semistable. Thus, there are Gieseker unstable instanton sheaves of
every rank.
\end{remark}

A little more can be said about rank 2 reflexive instanton sheaves on $\p3$ and
rank 4 locally-free instanton sheaves on $\p5$. 

\begin{proposition}\label{rk2p3}
Every rank 2 reflexive instanton sheaf on $\p3$ is locally-free and stable.
\end{proposition}
\begin{proof}
Hartshorne has shown that if $E$ is a rank 2 reflexive sheaf on $\p3$
with $c_3(E)=0$, then $E$ is locally-free \cite{Ha}, thus every rank 2
reflexive instanton sheaf on $\p3$ is locally-free. By Proposition \ref{h0},
we have that $H^0(E)=0$, hence $E$ is stable by Proposition \ref{m+rk2+rk3}.
\end{proof}

This result is sharp, in the sense that there are properly semistable
rank 2 torsion-free sheaves on $\p3$ and properly semistable rank 3
properly reflexive instanton sheaves on $\p3$:

\begin{example}\rm\label{dnoi}
Consider the monad:
\begin{equation} \label{ex-tf}
\op3(-1) \stackrel{\alpha}{\rightarrow} \op3^{\oplus4}
\stackrel{\beta}{\rightarrow} \op3(1)
\end{equation}
$$ \alpha = \left(\begin{array}{c} x_1 \\ x_2 \\ 0 \\ 0 \end{array}\right) 
~~{\rm and}~~ \beta= (-x_2 ~~ x_1 ~~ x_3 ~~ x_4) ~~. $$
Since $\alpha$ is injective provided $x_1,x_2\neq0$, its cohomology is a
rank 2 properly torsion-free instanton sheaf of charge 1. Moreover, $E$ is
not stable because it is a non-locally-free nullcorrelation sheaf
\cite[remark 1.2.1]{E}.

Finally, note that $E^*$ is a properly semistable rank 2 properly reflexive
sheaf on $\p3$ with $c_1(E^*)=0$; by Proposition \ref{rk2p3}, $E^*$ cannot be
instanton.
\end{example}

\begin{example}\rm
Set $w=5$ and $v=u=1$ and consider the monad:
\begin{equation} \label{ex-ref}
\op3(-1) \stackrel{\alpha}{\rightarrow} \op3^{\oplus5}
\stackrel{\beta}{\rightarrow} \op3(1)
\end{equation}
$$ \alpha = \left(\begin{array}{c} x_1 \\ x_2 \\ 0 \\ 0 \\ x_3 \end{array}\right) 
~~{\rm and}~~
\beta= (-x_2 ~~ x_1 ~~ x_3 ~~ x_4 ~~ 0) ~~.$$
It is easy to see that $\beta$ is surjective for all $[x_1:\cdots:x_4]\in\p3$,
while $\alpha$ is injective provided $x_1,x_2,x_3\neq0$. It follows that $E$ is
reflexive, but not locally-free; its singularity set is just the point $[0:0:0:1]\in\p3$.

In summary, $E$ is a rank 3 properly reflexive instanton sheaf of charge 1
on $\p3$. Note that $E$ is properly semistable, by Theorem \ref{highrank} and
Proposition \ref{nost}.
\end{example}

\begin{proposition}
Every rank 4 locally-free instanton sheaf $E$ on $\p5$ is stable.
\end{proposition}
\begin{proof}
Noting that $H^0(E)=0$ by Proposition \ref{h0}, the claim follows from
\cite[Theorem 3.6]{AO}. 
\end{proof}

Again, this result is sharp, in the sense that there exists a properly semistable
rank 4 properly reflexive instanton sheaf on $\p5$; it is also not true that every
rank $2n$ reflexive instanton sheaf on $\mathbb{P}^{2n+1}$ is locally-free or stable,
as Hartshorne's result could suggest.

\begin{example} \label{p5}\rm
Consider the cohomology $E$ of the monad:
$$ 0\to\op5(-1)\stackrel{\alpha}{\to}\op5^{\oplus 6}
\stackrel{\beta}{\to}\op5(1)\to0 $$
with the maps $\alpha$ and $\beta$ given by:
$$ \beta = \left( \begin{array}{cccccc} x_1~&x_2~&x_3~&x_4~&x_5~&x_6 \end{array} \right) ~~~~
\alpha = \left( \begin{array}{ccc} -x_2\\x_1\\-x_4\\x_3\\0\\0 \end{array} \right) .$$
Its degeneration locus is the line $\{x_1=\cdots=x_4=0\}$, so its cohomology is indeed
properly reflexive.

Finally, $E$ is properly semistable because it is a non-locally-free
nullcorrelation sheaf \cite[Remark 1.2.1]{E}, and $H^0(E^*)\neq0$.
\end{example}

\begin{remark}\rm
It seems reasonable to conjecture that every rank $2n$ locally-free instanton sheaf on
$\mathbb{P}^{2n+1}$ is stable. In support of this conjecture, see Lemma \ref{simp}
and Proposition \ref{k1st} below. Results in this direction were also obtained in
\cite{AO} for {\em symplectic} mathematical instanton bundles.
\end{remark}

\begin{example}\label{ex4}\rm
Indecomposable rank 2 locally-free sheaves on $\pn$, $n\geq4$, have been
extremely difficult to construct, and linear monads do not help with this
problem. However, stable rank 2 reflexive sheaves on $\pn$ are easy to
construct. Indeed, Fl\o ystad's theorem guarantees the existence of linear
monad:
$$ 0 \to \opn(-1)^{\oplus n+a-3} \stackrel{\alpha}{\to} \opn^{\oplus n+2a-1}
\stackrel{\beta}{\to} \opn(1)^{\oplus a} \to 0 ~~ (a\ge1) ~~$$
whose cohomology is a rank $2$ reflexive linear sheaf with $c_1(E)=n-3$.
To see that it is also stable, note that $E_\eta=E(k)$ for some $k\leq-1$,
thus $H^0(E_\eta)=0$ and it follows from Lemma \ref{rk2tf} that $E$ must be
stable.

It is interesting to contrast the existence of such stable rank 2 reflexive
sheaves on $\pn$ with Hartshorne's conjecture: there are no indecomposable
rank 2 locally-free sheaves on $\pn$ for $n\ge 7$ \cite{Ha1}.

Furthermore, this example implies that Kumar, Peterson and Rao's result (Theorem \ref{mkr})
is sharp, in the sense it cannot be extended to more general sheaves: there are
rank 2 reflexive sheaves $E$ on $\pn$ with $H^p_*(E)=0$ for $2\leq p\leq n-2$
which do not split as a sum of rank 1 sheaves.

For instance, with $n=4$ and $a=1$, we get the monad:
$$ 0\to\op4(-1)^{\oplus 2}\stackrel{\alpha}{\to}\op4^{\oplus 5}
\stackrel{\beta}{\to}\op4(1)\to0 $$
with the maps $\alpha$ and $\beta$ given by:
$$ \beta = \left( \begin{array}{ccccc} x_1~&x_2~&x_3~&x_4~&x_5 \end{array} \right) ~~~~
\alpha = \left( \begin{array}{cc} -x_2&-x_5\\x_1&x_3\\-x_4&-x_2\\x_3&0\\0&x_1 \end{array} \right) ,$$
where $[x_1:\dots:x_5]$ are homogeneous coordinates in $\p4$.
Note that the degeneration locus of this monad is given by the union of two lines:
$$ \Sigma(E)=\{x_1=x_2=x_3=x_4=0\}\cup\{x_1=x_2=x_3=x_5=0\}~. $$
\end{example}

One can thus hope to construct stable rank 2 locally-free sheaves on $\pn$ via
some mechanism that turns reflexive into locally-free sheaves without introducing
new global sections.


\paragraph{Lifting of instantons.}
It is known that for every locally-free instanton sheaf $E$ of trivial
splitting type on $\p2$, there exists a locally-free instanton sheaf
$\tilde{E}$ of trivial splitting type on $\p3$ and a hyperplane $\wp\subset\p3$
such that $\tilde{E}|_{\wp}\simeq E$ \cite{D}. $\tilde{E}$ is called a {\em lifting}
of $E$.

It would be interesting to see whether this generalizes to higher
dimensional projective spaces and/or to more general sheaves; more precisely,
we propose the following conjecture:

\begin{conjecture}
If $E$ is a rank $r\ge2n$ locally-free instanton sheaf of trivial splitting type
on $\mathbb{P}^{2n}$ of charge $c$, then there is a rank $r$ locally-free instanton
sheaf $\tilde{E}$ on $\mathbb{P}^{2n+1}$ of charge $c$ and a hyperplane
$\wp\subset\mathbb{P}^{2n+1}$ such that $\tilde{E}|_{\wp}\simeq E$.
\end{conjecture}

Donaldson's argument in \cite{D} for the $n=1$ case is ``unashamedly computational",
relying in the correspondence between instanton sheaves of trivial splitting type on
$\p2$ and $\p3$ and solutions of the ADHM equations, and it is not clear how to generalize
it for $n\ge2$. As far as we know, there is no alternative, more conceptual proof of
Donaldson's result.

It is not difficult to see that if a locally-free instanton sheaf on $\mathbb{P}^{2n}$ can
be lifted to $\mathbb{P}^{2n+1}$, then the lifted sheaf must be instanton. The hard part
is determining whether a given instanton given sheaf can be lifted;
the condition that the sheaf on $\mathbb{P}^{2n}$ is of trivial splitting type might
be crucial here.

\begin{proposition}
Let $E$ be a locally-free instanton sheaf on $\mathbb{P}^{n}$ ($n\ge4$) of
charge $c$. If there is a locally-free sheaf $\tilde{E}$ on $\mathbb{P}^{n+1}$,
and a hyperplane $\wp\subset\mathbb{P}^{n+1}$ such that $\tilde{E}|_{\wp}\simeq E$, 
then $\tilde{E}$ is also an instanton sheaf of charge $c$.
\end{proposition}
\begin{proof}
The desired result follows from the restriction sequence:
$$ 0 \to \tilde{E}(k-1) \to \tilde{E}(k) \to E(k) \to 0 ~~. $$
together with repeated use of Serre's duality and Serre's vanishing
theorem.

Since $H^0(E(k))=0$ for all $k\le-1$, we get that 
$$ H^0(\tilde{E}(k-1)) \simto H^0(\tilde{E}(k)) $$
for all $k\le-1$. Hence also $H^0(\tilde{E}(-1))=0$.

Similarly, since $H^n(E(k))=0$ for all $k\ge-n$, we get that
$$ H^{n+1}(\tilde{E}(k-1)) \simto H^{n+1}(\tilde{E}(k)) $$
for all $k\ge-n$. Hence also $H^{n+1}(\tilde{E}(-n-1))=0$.

Since $H^1(E(k))=0$ for all $k\le-2$, we get that
$$ H^1(\tilde{E}(k-1)) \to H^1(\tilde{E}(k))\to0 $$
for all $k\le-2$. But $H^1(\tilde{E}(l))=0$ for $l\ll0$ by Serre's
vanishing theorem, thus it follows that $H^1(\tilde{E}(-2))=0$.

Similarly, since $H^{n-1}(E(k))=0$ for all $k\ge1-n$, we get that
$$ 0\to H^{n}(\tilde{E}(k-1)) \to H^{n}(\tilde{E}(k)) $$
for all $k\ge1-n$. But $H^1(\tilde{E}(l))=0$ for $l\gg0$ by Serre's
vanishing theorem, thus it follows that $H^n(\tilde{E}(-n))=0$.
This completes the proof for the case $n=2$.

Since $n\ge4$, we have that $H^p(E(k))=0$ for $2\leq p\leq n-2$ and all $k$;
thus
$$ H^2(\tilde{E}(k-1)) \simto H^2(\tilde{E}(k)) $$
for all $k\le-2$. Again Serre's vanishing theorem forces $H^2(\tilde{E}(k))=0$
for all $k\le-2$. Moreover, we have that
$$ H^2(\tilde{E}(k-1)) \to H^2(\tilde{E}(k)) \to 0 $$
for all $k$, hence $H^2(\tilde{E}(k))=0$ for all $k$.

Finally, for $n\ge5$ we have that
$$ H^p(\tilde{E}(k-1)) \simto H^p(\tilde{E}(k)) $$
for $3\le p\le n-2$ and all $k$. It follows that 
$H^p(\tilde{E}(k))=0$ for $3\le p\le n-2$ and all $k$,
which completes the proof of the first statement.

Setting $k=-1$ on the restriction sequence, we get
$h^1(\tilde{E}(-1))=h^1(E(-1))$, showing that the charge is preserved.
\end{proof}


\section{Simplicity of linear sheaves}\label{simple}
Recall that a torsion-free sheaf $E$ on a projective variety $X$ is said to be
{\em simple} if $\dim{\rm Ext}^0(E,E)=1$. Every stable torsion-free sheaf on
$\pn$ is simple. 

\begin{theorem}\label{simple.thm}
If $E$ is the cohomology of the linear monad
$$ 0\to  V\otimes\opn(-1) \stackrel{\alpha}{\longrightarrow}
W\otimes\opn \stackrel{\beta}{\longrightarrow} U\otimes\opn(1) \to 0 ~, $$
for which $K=\ker\beta$ is simple, then $E$ is simple.
\end{theorem}
\begin{proof}
Applying ${\rm Ext}^*(\cdot,E)$ to the sequence
$$ 0 \to V\otimes\opn(-1) \to K \to E \to 0 ~~,$$
we get
\begin{equation}\label{ee}
0\to{\rm Ext}^0(E,E)\to{\rm Ext}^0(K,E)\to\cdots ~~.
\end{equation}
Now applying ${\rm Ext}^*(K,\cdot)$ we get:
$$ V\otimes{\rm Ext}^0(K,\opn(-1)) \to {\rm Ext}^0(K,K)\to
{\rm Ext}^0(K,E) \to V\otimes{\rm Ext}^1(K,\opn(-1)). $$
But it follows from the dual of sequence (\ref{ker1}) with $k=1$ that
$h^0(K^*(-1))=h^1(K^*(-1))=0$, thus
$$ \dim{\rm Ext}^0(K,E)=\dim{\rm Ext}^0(K,K)=1 $$
because $K$ is simple. It then follows from (\ref{ee}) that
$E$ is also simple.
\end{proof}

As a consequence of \cite[Theorem 2.8(a)]{AO}, we have
in particular the following generalization of \cite[Theorem 2.8(b)]{AO}:

\begin{lemma}\label{simp}
Every rank $n-1$ linear sheaf on $\pn$ is simple. 
\end{lemma}

The above result is sharp, in the sense that there are rank $n$ instanton
sheaves on $\pn$ which are not simple. For example, recall that a rank 2
locally-free sheaf is simple if and only if it is stable; since every 
rank 2 instanton sheaf on $\p2$ of charge 1 is properly semistable
(by Proposition \ref{nost}), it follows that these are not simple, as
desired.


\section{Moduli spaces of instanton sheaves}\label{moduli}

Let ${\cal I}_{\mathbb{P}^n}(r,c)$ denote the moduli space
of equivalence classes of rank $r$ instanton sheaves of charge
$c$ on $\pn$. Let ${\cal I}_{\mathbb{P}^n}^{\rm lf}(r,c)$ denote
the open subset of ${\cal I}_{\mathbb{P}^n}(r,c)$ consisting
of locally-free sheaves. Note that ${\cal I}_{\mathbb{P}^n}^{\rm lf}(r,c)$
might be empty even though ${\cal I}_{\mathbb{P}^n}(r,c)$ is not.

Very little is known in general about ${\cal I}_{\mathbb{P}^n}(r,c)$;
research so far has concentrated on ${\cal I}_{\mathbb{P}^{2n+1}}^{\rm lf}(2n,c)$.
Here is a summary of some of the known facts:\label{lit}
\begin{itemize}
\item ${\cal I}_{\mathbb{P}^{2n+1}}^{\rm lf}(2n,c)$ is affine \cite{CO};
\item ${\cal I}_{\mathbb{P}^{2n+1}}^{\rm lf}(2n,1)$ is an open subset of
$\mathbb{P}^{n(2n+1)-1}$ \cite{E};
\item ${\cal I}_{\mathbb{P}^{2n+1}}^{\rm lf}(2n,2)$ is an irreducible,
smooth variety of dimension $4n^2+12n-3$ \cite{AO};
\item ${\cal I}_{\mathbb{P}^{3}}^{\rm lf}(2,c)$ is an irreducible,
smooth variety of dimension $8c-3$ for $1\leq c\leq5$, see \cite{CTT}
and the references therein;
\item ${\cal I}_{\mathbb{P}^{2n+1}}^{\rm lf}(2n,c)$ is singular for
all $n\geq2$ and $c\geq3$ \cite{MROL}.
\end{itemize}
The smoothness of ${\cal I}_{\mathbb{P}^{3}}^{\rm lf}(2,c)$ for arbitrary
charge $c$ is still an open problem; it is known however that its
closure in the moduli space of semistable locally-free sheaves with Chern
character $2-cH^2$ is in general singular, see \cite{AO2}.

In this section, we will generalize the second statement, study
the moduli spaces of instanton sheaves on $\p2$ and conclude with a
general conjecture that generalizes the first statement.


\paragraph{Instanton sheaves and nullcorrelation sheaves}\label{nc}

Recall that a nullcorrelation sheaf $N$ on $\pn$ is a rank $n-1$
torsion-free sheaf defined by the short exact sequence:
$$ 0\to\opn(-1)\stackrel{\sigma}{\to}\Omega^1_{\pn}(1)\to N \to 0 ~~,$$
where $\sigma\in H^0(\Omega^1_{\pn}(2))=\Lambda^2H^0(\opn(1))$. If
$n$ is odd and $\sigma$ is generic, then $N$ is locally-free. If $n$
is even, then $N$ is never locally-free; however, the generic one is
reflexive if $n\geq4$.

It is easy to see that any nullcorrelation sheaf is instanton of charge 1.
The converse is also true: every rank $n-1$ torsion-free instanton sheaf
of charge 1 is a nullcorrelation sheaf. Indeed, let $E$ be a rank $n-1$
instanton sheaf of charge 1, so that it is the cohomology of the
sequence:
$$ 0\to\opn(-1)\stackrel{\alpha}\to\opn^{\oplus(n+1)}
\stackrel{\beta}\to\opn(1)\to 0 ~~. $$
Comparing this with the Euler sequence, it follows that $\ker\beta$ coincides
$\Omega^1_{\pn}(1)$, up to an automorphism of $\pn$. Thus $E$ fits into the
sequence:
$$ 0\to\opn(-1)\stackrel{\alpha}{\to}\Omega^1_{\pn}(1)\to E \to 0 ~~,$$
hence $E$ is nullcorrelation.

It follows from the above correspondence and the fact that every nullcorrelation
sheaf is simple that any nullcorrelation sheaf is completely determined by a section
$\sigma\in H^0(\Omega^1_{\pn}(2))$. Hence, we have that the set of equivalence
classes of nullcorrelation sheaves on $\pn$ is exactly
$\mathbb{P}\left(H^0(\Omega^1_{\pn}(2))\right)$. Thus we conclude:

\begin{theorem}
${\cal I}_{\mathbb{P}^{n}}(n-1,1)\simeq
\mathbb{P}^{\frac{n(n+1)}{2}-1}$.
\end{theorem}

It is known that every nullcorrelation locally-free sheaf is stable,
while nullcorrelation sheaves which are not locally-free are not stable
\cite[Remark 1.2.1]{E}; they are however properly semistable. In particular,
we have:

\begin{proposition}\label{k1st}
Every rank $2n$ locally-free instanton sheaf of charge 1 on $\mathbb{P}^{2n+1}$
is stable.
\end{proposition}


\paragraph{Moduli spaces of instanton sheaves on $\p2$.}

The simplest possible instanton sheaves are the rank 1 instanton
sheaves on $\p2$. It is not difficult to see that such sheaves are
exactly the ideals of points in $\p2$.

Indeed, let $Z$ be a closed zero dimensional subscheme in $\p2$, and
let $I_Z$ denote its ideal sheaf; it fits into the sequence:
\begin{equation}\label{ideal}
0 \to I_Z \to \op2 \to {\cal O}_Z \to 0 ~~.
\end{equation}
After tensoring with $\op2(k)$, it follows
that $H^0(\p2,I_Z(k))=0$ for $k\leq-1$, and $H^2(\p2,I_Z(k))=0$ for
$k\geq-2$, so $I_Z$ is indeed instanton. Moreover, the charge of $I_Z$
is just the length of $Z$.

Conversely, let $E$ be a rank 1 instanton sheaf of charge $c$ on $\p2$.
Then $E^{**}$ is a rank 1 locally-free sheaf with $c_1(E^{**})=c_1(E)=0$,
so $E^{**}=\op2$. Thus $E$ is the ideal sheaf associated with the
zero-dimensional scheme $\op2/E$, whose length is equal to the charge of
$E$.

In other words, there is a 1-1 correspondence between rank 1 instanton sheaves
of charge $c$ on $\p2$ and closed zero dimensional subschemes of length $c$ in
$\p2$. This gives us the following identity:

\begin{theorem}
${\cal I}_{\mathbb{P}^{2}}(1,c) \simeq (\p2)^{[c]}$.
\end{theorem}

\begin{remark}\rm
Let $E$ be a locally-free instanton sheaf of charge $c$; by tensoring sequence
(\ref{ideal}) with $E(k)$, it is easy to see that $E_Z=E\otimes I_Z$ is a properly
torsion-free instanton sheaf of charge $c+r\cdot~{\rm lenght}(Z)$. Moreover, if
$H^0(\p2,E)=0$, then $H^0(\p2,E_Z)=0$. Is it true that every properly torsion-free
instanton is a locally-free instanton sheaf tensored by an ideal of points?
\end{remark}

As we have seen in Section \ref{ss}, rank $2$ and $3$ instanton sheaves on
$\p2$ are in 1-1 correspondence with semistable torsion-free sheaves with
zero first Chern class. So we have:

\begin{corollary}
For $r=2,3$, 
${\cal I}_{\mathbb{P}^{2}}(r,c) \simeq {\cal M}_{\mathbb{P}^{2}}(r,0,c)$,
the moduli space of rank $r$ semistable torsion-free sheaves $E$ on $\p2$ with
$c_1(E)=0$ and $c_2(E)=c$.
In particular, ${\cal I}_{\mathbb{P}^{2}}(2,c)$ and ${\cal I}_{\mathbb{P}^{2}}(3,c)$
are quasi-projective varieties.
\end{corollary}

In general, we conjecture that ${\cal I}_{\pn}(r,c)$ is always a quasi-projective
variety.


 \end{document}